\documentclass[11pt]{amsart}
\usepackage{amsmath,amssymb,amsfonts,amsthm,latexsym}

\newcommand{\hh}{\ensuremath{\mathbf{h}}}

\newcommand{\pfrac}[2]{{\left(\frac{#1}{#2}\right)}}

\newcommand{\be}{\begin{equation}}
\newcommand{\ee}{\end{equation}}
\newcommand{\benn}{\begin{equation*}}   
\newcommand{\eenn}{\end{equation*}}

\renewcommand{\(}{\left(}
\renewcommand{\)}{\right)}

\begin{document}

\title{Simple proof of Gallagher's singular series sum estimate}

\author{Kevin Ford}


\date{October 13, 2007}

%

\maketitle

P. X. Gallagher (Mathematika {\bf 23} (1976), 4--9) 
proved an estimate for the average
of the singular series associated with $r$-tuples of linear forms, namely
\be\label{G}
\sum_{\substack{0\le h_1, \ldots, h_r \le h \\ h_1,\ldots, h_r
    \text{ distinct}}} \prod_{p} \( 1 -
\frac{\nu_p(\hh)}{p}\)\(1-\frac{1}{p}\)^{-r} \sim h^r
\ee
for each fixed $r$, where $\nu_p(\hh)$ is the number of residue
classes modulo $p$ occupied by the numbers $h_1,\ldots, h_r$.  We give
a simpler proof of this result below, with a worse error
estimate than Gallagher obtained.  The constants in all $O$-terms
may depend on $r$.

Put $y=\frac12 \log h$.
We first note that $\nu_p(\hh)=r$ if $p\nmid H$, where $H=\prod_{i<j}
|h_i-h_j|$.  The number of prime factors of $H$ is $O(\log H/\log\log H)=O(\log
h/\log\log h)$.   For any $h_1,\ldots,h_r$, we therefore have
\begin{align*}
\prod_{p>y} \( 1 - \frac{\nu_p(\hh)}{p}\)\(1-\frac{1}{p}\)^{-r} &=
\prod_{p|H, p>y} \(1+O\pfrac{1}{p}\) \prod_{p\nmid H,p>y} \(1+O\pfrac{1}{p^2}\) \\
&= 1 + O\pfrac{\log h}{y\log\log h} = 1 + O\pfrac{1}{\log\log h}.
\end{align*}
Thus, the left side of \eqref{G} is equal to $AB$, where
$$
A=\(1 + O\pfrac{1}{\log\log h}\) \prod_{p\le y}
\(1-\frac{1}{p}\)^{-r}, \quad B=
\sum_{\substack{0\le h_1, \ldots, h_r \le h \\ h_1,\ldots, h_r
    \text{ distinct}}} \prod_{p\le y} \(1-\frac{\nu_p(\hh)}{p}\).
$$
We have $B=O(h^{r-1}) + B'$, where $B'$ is the corresponding sum
without the condition that $h_1,\ldots,h_r$ are distinct.
Let $P=\prod_{p\le y} p$ and note that $P=e^{y+o(y)}=h^{1/2+o(1)}$.
The product in $B$ is $1/P$ times the
number of $n$, $0\le n<P$, satisfying $(\prod_i (n+h_i),P)=1$.
Threrefore,
\begin{align*}
B' &= \sum_{0\le h_1,\ldots,h_r\le h} \frac{1}{P} \sum_{n=0}^{P-1}
\prod_{i=1}^r \sum_{d_i|(n+h_i,P)} \mu(d_i) \\
&=\frac{1}{P} \sum_{n=0}^{P-1} \sum_{d_1,\ldots,d_r|P} \mu(d_1)\cdots
\mu(d_r) \prod_{i=1}^r \( \frac{h}{d_i}+O(1)\) \\
&=h^r  \sum_{d_1,\ldots,d_r|P} \frac{\mu(d_1) \cdots
  \mu(d_r)}{d_1\cdots d_r} + O\( h^{r-1}   \sum_{d_1,\ldots,d_{r-1}|P} 
\frac{1}{d_1\cdots d_r} \) \\
&=h^r \prod_{p\le y} \(1-\frac{1}{p}\)^r + O(h^{r-1+o(1)}).
\end{align*}
Combined with the expression for $A$, this proves \eqref{G}.

\end{document}